\documentclass[12pt]{article}

\usepackage[letterpaper,margin=1in]{geometry}
\usepackage{amsmath,amssymb,amsthm,mathtools}
\usepackage{array}
\usepackage{booktabs}
\usepackage{float}
\usepackage{microtype}
\usepackage{mathrsfs}
\usepackage{xcolor}
\usepackage{hyperref}

\hypersetup{
  colorlinks=true,
  linkcolor=blue!55!black,
  citecolor=blue!55!black,
  urlcolor=blue!55!black
}
\numberwithin{equation}{section}
\allowdisplaybreaks

\newtheorem{theorem}{Theorem}[section]

\newtheorem{lemma}[theorem]{Lemma}

\newcommand{\dd}{\,\mathrm d}
\newcommand{\Eop}{\mathcal E}
\newcommand{\Dop}{\mathscr D}

\title{\textbf{Even and Odd Integer Values of the Secondary Zeta Function}}
\author{Artur Kawalec}
\date{}

\begin{document}
\maketitle

\begin{abstract}
The secondary zeta function is defined by
\[
Z(s)=\sum_{\gamma>0}\frac{1}{\gamma^{s}},
\]
where $\tfrac12+i\gamma$ runs over the non-trivial zeros of the Riemann zeta function with positive ordinate, and the Riemann hypothesis is assumed throughout. In this article, we derive convergent integral representations of $Z(s)$ at positive even and odd integer argument, expressed through a differential operator acting on the logarithm of the Riemann $\xi$-function. Define
\[
L(x)=\log\!\left(\frac{\xi(x+\tfrac12)}{\xi(\tfrac12)}\right),
\qquad
\Eop=x\frac{\dd}{\dd x},
\qquad
\Dop_k=\Eop\prod_{j=1}^{k}(\Eop-2j).
\]
Then, for $k\geq1$,
\[
Z(2k)=\frac{(-1)^k}{2^k(k-1)!}
\int_{0}^{\infty}\frac{\Dop_k L(x)}{x^{2k+1}}\,\dd x
\]
and
\[
Z(2k+1)=\frac{(-1)^k}{\pi\,(2k-1)!!}
\int_{0}^{\infty}\frac{\Dop_k L(x)}{x^{2k+2}}\,\dd x.
\]
In the last part, we write these formulas out explicitly for $Z(k)$, $2\leq k\leq 7$, both in terms of $L(x)$ and in terms of $\xi(x+\tfrac12)$ for reference. We also derive a limit formula for $\gamma_1$, the ordinate of the first nontrivial zero, and demonstrate its convergence numerically.
\end{abstract}

\section{Introduction}

The secondary zeta function associated with the ordinates of the non-trivial zeros of the Riemann zeta function has a long history; see Delsarte \cite{Delsarte1966}, Chakravarty \cite{Chakravarty1970}, Ivi\'c \cite{Ivic2003}, Voros \cite{Voros2003}, and Arias de Reyna \cite{Arias2020}.  We assume the Riemann hypothesis (RH) throughout and write
\[
\rho=\frac12+i\gamma,
\qquad \gamma>0,
\]
so that
\begin{equation}
Z(s)=\sum_{\gamma>0}\frac{1}{\gamma^s},
\qquad \Re s>1.
\label{eq:defZ}
\end{equation}
The function $Z(s)$ extends meromorphically to the complex plane; its analytic continuation and special values have been studied in, among other sources, \cite{Voros2003,BondarenkoIvic2018,Arias2020,Kawalec2024}.

A recent historical note by Arias de Reyna revisits Riemann's calculations of logarithmic derivatives at the central point, including the second derivative underlying the evaluation of $Z(2)$ by Riemann \cite{AriasRiemann2026}.

Table~\ref{tab:bench} gives numerical reference values for all integer arguments from $2$ through $11$.  The even values were checked independently from the centered $\xi$-derivative formula in Section~\ref{sec:even}, while all entries can also be obtained from the analytic-continuation algorithm of Arias de Reyna \cite{Arias2020}.

\begin{table}[htbp]
\centering
\caption{Reference values of the secondary zeta function to $50$ decimal digits.}
\label{tab:bench}
\begingroup
\small
\setlength{\tabcolsep}{5pt}
\renewcommand{\arraystretch}{1.13}
\begin{tabular}{@{}c>{\ttfamily}l@{}}
\toprule
$ n $ & \multicolumn{1}{c}{$Z(n)$} \\
\midrule
2  & 0.02310499311541897078893381043033901400338176039742  \\
3  & 0.00072954827270970421587551856909397050335150570355  \\
4  & 0.00003717259928526968616486626247174057845365088973  \\
5  & 0.00000223118869950210332864062869183719337607643108  \\
6  & 0.00000014417393140097327969538155609482090703688300  \\
7  & 0.00000000967534454270235040871996562746154140230137  \\
8  & 0.00000000066303168025299086987327208196135724847369  \\
9  & 0.00000000004599191239289486296927465194968984893450  \\
10 & 0.00000000000321366415061660121610211659983465514156  \\
11 & 0.00000000000022556506251559664289836745418802706353  \\
\bottomrule
\end{tabular}
\endgroup
\end{table}

The central observation is that the centered Hadamard product of the Riemann $\xi$-function contains the even secondary-zeta values as Taylor coefficients, while a regularized Mellin transform of the same logarithm recovers the odd values.  An Euler differential operator can moreover be chosen to annihilate exactly the low-order Taylor terms that otherwise have to be subtracted explicitly.

\section{Even values from the Riemann \texorpdfstring{$\xi$}{xi}-function}
\label{sec:even}

Let $\xi(s)$ denote the Riemann $\xi$-function in its standard normalization, so that $\xi(s)=\xi(1-s)$.  Under RH, pairing the zeros with ordinates $\pm\gamma$ gives the centered product
\begin{equation}
\frac{\xi(x+\tfrac12)}{\xi(\tfrac12)}
=\prod_{\gamma>0}
\left(1+\frac{x^2}{\gamma^2}\right).
\label{eq:centeredproduct}
\end{equation}
For real $x$, interpreted as displacement from the central point $s=\tfrac12$, define
\begin{equation}
\boxed{
L(x)=\log\!\left(\frac{\xi(x+\tfrac12)}{\xi(\tfrac12)}\right).
}
\label{eq:Ldef}
\end{equation}
It follows that
\begin{equation}
L(x)
=\sum_{\gamma>0}\log\!\left(1+\frac{x^2}{\gamma^2}\right)
=\sum_{m=1}^{\infty}(-1)^{m+1}\frac{Z(2m)}{m}x^{2m},
\label{eq:Lseries}
\end{equation}
for $|x|<\gamma_1$.  Comparing Taylor coefficients gives
\begin{equation}
\boxed{
Z(2m)=\frac{(-1)^{m+1}}{2(2m-1)!}\,L^{(2m)}(0)
=\frac{(-1)^{m+1}}{2(2m-1)!}
\left.\frac{\dd^{2m}}{\dd s^{2m}}\log\xi(s)\right|_{s=1/2}.
}
\label{eq:Zevenxi}
\end{equation}
For reference, the equivalent functional-equation form used by Voros may be written as follows.  Set
\begin{equation}
D_{2m}:=
\left.
\frac{\dd^{2m}}{\dd s^{2m}}\log\zeta(s)
\right|_{s=1/2}.
\label{eq:D2m}
\end{equation}
Then
\begin{equation}
\boxed{
Z(2m)=(-1)^m\left[
-\frac{D_{2m}}{2(2m-1)!}
-\frac14\Bigl((2^{2m}-1)\zeta(2m)+2^{2m}\beta(2m)\Bigr)
+2^{2m}
\right],
}
\label{eq:VorosEvenFunctional}
\end{equation}
where
\(
\beta(s)=\sum_{n=0}^{\infty}(-1)^n(2n+1)^{-s}
\)
is the Dirichlet beta function.  Formula \eqref{eq:VorosEvenFunctional} results from expanding the elementary and gamma factors entering the functional equation of $\zeta(s)$; see Voros \cite{Voros2003}.  Thus the same single function $L(x)$ supplies all of the even coefficients used in the odd-value formulas below.

\section{A convergent integral formula for the even values}
\label{sec:evenintegral}

Formula~\eqref{eq:Zevenxi} recovers $Z(2k)$ locally from the Taylor
coefficient of $L(x)$ at the origin.  There is also a convergent integral
counterpart that is directly parallel to the differential-operator formula
for the odd values.  After the lower Taylor powers are removed from the
Mellin integrand with weight $x^{-2k-1}$, the surviving $x^{2k}$ term still
produces a logarithmic singularity at $x=0$.  Applying an Euler operator
that also annihilates this borderline term yields an ordinary convergent
integral.

For this purpose, set
\[
\Eop=x\frac{\dd}{\dd x},
\qquad
\Dop_k=\Eop\prod_{j=1}^{k}(\Eop-2j).
\]
Thus $\Dop_k$ annihilates every monomial $x^{2m}$ with $1\leq m\leq k$.

\begin{theorem}[Even-value differential-operator integral]
For every integer $k\geq1$,
\begin{equation}
\boxed{
Z(2k)=\frac{(-1)^k}{2^k(k-1)!}
\int_{0}^{\infty}
\frac{\Dop_kL(x)}{x^{2k+1}}\,\dd x.
}
\label{eq:evenoperatorintegral}
\end{equation}
Equivalently, using the definition of $L(x)$ explicitly,
\begin{equation}
\boxed{
Z(2k)=\frac{(-1)^k}{2^k(k-1)!}
\int_{0}^{\infty}
\frac{\Dop_k\log\!\left(\dfrac{\xi(x+\tfrac12)}{\xi(\tfrac12)}\right)}
{x^{2k+1}}\,\dd x.
}
\label{eq:evenoperatorxi}
\end{equation}
\end{theorem}

\begin{proof}
For a single positive ordinate $\gamma$, write
\[
f_\gamma(x)=\log\!\left(1+\frac{x^2}{\gamma^2}\right).
\]
With $t=x^2/\gamma^2$, one has $\Eop=2t\,\dd/\dd t$.  A direct induction
therefore gives
\begin{equation}
\Dop_k f_\gamma(x)
=(-1)^k2^{k+1}k!\,
\frac{x^{2k+2}}{(\gamma^2+x^2)^{k+1}}.
\label{eq:evenkernel}
\end{equation}
Consequently,
\[
(-1)^k\frac{\Dop_k f_\gamma(x)}{x^{2k+1}}
=2^{k+1}k!\frac{x}{(\gamma^2+x^2)^{k+1}},
\]
which is non-negative.  Since
\(
L(x)=\sum_{\gamma>0}f_\gamma(x)
\), Tonelli's theorem permits the sum over $\gamma$ to be interchanged
with the integral, and
\[
\int_0^\infty\frac{x}{(\gamma^2+x^2)^{k+1}}\,\dd x
=\frac{1}{2k\,\gamma^{2k}}.
\]
It follows that
\[
\int_{0}^{\infty}
\frac{\Dop_kL(x)}{x^{2k+1}}\,\dd x
=(-1)^k2^k(k-1)!\,Z(2k),
\]
which proves \eqref{eq:evenoperatorintegral}.  The same calculation also
proves absolute convergence.  Near $x=0$, the operator $\Dop_k$ removes
the terms $x^2,x^4,\ldots,x^{2k}$, so the resulting integrand is $O(x)$;
at the upper endpoint, Tonelli's identity gives a finite total integral
because $Z(2k)<\infty$.
\end{proof}

The first three instances are
\begin{align}
Z(2)
&=\frac12\int_0^\infty
\left(\frac{L'(x)}{x^2}-\frac{L''(x)}{x}\right)\dd x,
\label{eq:Z2evenintegral}\\
Z(4)
&=\frac14\int_0^\infty
\left(
\frac{L'''(x)}{x^2}
-\frac{3L''(x)}{x^3}
+\frac{3L'(x)}{x^4}
\right)\dd x,
\label{eq:Z4evenintegral}\\
Z(6)
&=-\frac1{16}\int_0^\infty
\left(
\frac{L^{(4)}(x)}{x^3}
-\frac{6L'''(x)}{x^4}
+\frac{15L''(x)}{x^5}
-\frac{15L'(x)}{x^6}
\right)\dd x.
\label{eq:Z6evenintegral}
\end{align}
As with the odd-value formulas below, each displayed parenthesized
combination must be kept intact: the apparent singularities of its
individual terms cancel at the lower endpoint.

\section{A Taylor-remainder integral for every odd value}
\label{sec:remainder}

\begin{theorem}[Odd-value remainder integral]
For every integer $k\geq1$,
\begin{equation}
\boxed{
Z(2k+1)=\frac{(-1)^k(2k+1)}{\pi}
\int_{0}^{\infty}
\frac{L(x)-P_k(x)}{x^{2k+2}}\,\dd x,
}
\label{eq:remainderintegral}
\end{equation}
where
\begin{equation}
P_k(x)=\sum_{m=1}^{k}(-1)^{m+1}\frac{Z(2m)}{m}x^{2m}.
\label{eq:Pk}
\end{equation}
\end{theorem}

\begin{proof}
For one ordinate $\gamma$, put $y=x/\gamma$ and define
\[
r_k(y)=\log(1+y^2)-\sum_{m=1}^{k}(-1)^{m+1}\frac{y^{2m}}{m}.
\]
A finite geometric sum gives
\[
r_k'(y)=\frac{2(-1)^k y^{2k+1}}{1+y^2}.
\]
The boundary terms vanish in an integration by parts, and hence
\[
\int_0^{\infty}\frac{r_k(y)}{y^{2k+2}}\,\dd y
=\frac{1}{2k+1}\int_0^{\infty}\frac{r_k'(y)}{y^{2k+1}}\,\dd y
=\frac{(-1)^k\pi}{2k+1}.
\]
Scaling by $x=\gamma y$ supplies the factor $\gamma^{-(2k+1)}$.  Moreover, $(-1)^kr_k(y)\geq0$ for $y\geq0$, so Tonelli's theorem justifies interchanging the sum over $\gamma$ with the integral.  Summing over all positive ordinates proves \eqref{eq:remainderintegral}.
\end{proof}

\section{An annihilating differential operator}
\label{sec:operator}

Introduce the centered Euler operator
\begin{equation}
\Eop=x\frac{\dd}{\dd x}
\label{eq:Euler}
\end{equation}
and, for $k\geq1$, define
\begin{equation}
\boxed{
\Dop_k=\Eop\prod_{j=1}^{k}(\Eop-2j).
}
\label{eq:Dop}
\end{equation}
Because
\[
\Eop x^{2m}=2m x^{2m},
\]
the operator $\Dop_k$ annihilates the first $k$ centered even monomials.

\begin{lemma}
For
\[
f_\gamma(x)=\log\!\left(1+\frac{x^2}{\gamma^2}\right),
\]
one has
\begin{equation}
\boxed{
\Dop_k f_\gamma(x)
=(-1)^k2^{k+1}k!\,
\frac{x^{2k+2}}
{(\gamma^2+x^2)^{k+1}}.
}
\label{eq:Dkernel}
\end{equation}
\end{lemma}

\begin{proof}
Set $t=x^2/\gamma^2$.  Then $\Eop=2t\,\dd/\dd t$ and
\[
\Eop f_\gamma=\frac{2t}{1+t}.
\]
An induction gives
\[
\prod_{j=1}^{k}\left(t\frac{\dd}{\dd t}-j\right)\frac{t}{1+t}
=(-1)^k k!\frac{t^{k+1}}{(1+t)^{k+1}},
\]
which is equivalent to \eqref{eq:Dkernel}.
\end{proof}

\begin{theorem}[Differential-operator integral]
For every integer $k\geq1$,
\begin{equation}
\boxed{
Z(2k+1)=\frac{(-1)^k}{\pi(2k-1)!!}
\int_{0}^{\infty}
\frac{\Dop_kL(x)}{x^{2k+2}}\,\dd x.
}
\label{eq:mainoperatorintegral}
\end{equation}
Here $(2k-1)!!=1\cdot3\cdots(2k-1)$.
\end{theorem}

\begin{proof}
By \eqref{eq:Dkernel},
\[
\frac{\Dop_kL(x)}{x^{2k+2}}
=(-1)^k2^{k+1}k!
\sum_{\gamma>0}
\frac{1}{(\gamma^2+x^2)^{k+1}}.
\]
After multiplication by $(-1)^k$, all terms are non-negative, so Tonelli's theorem applies.  The beta integral gives
\[
\int_{0}^{\infty}
\frac{\dd x}{(\gamma^2+x^2)^{k+1}}
=\frac{\sqrt\pi\,\Gamma(k+\tfrac12)}{2\Gamma(k+1)}\gamma^{-(2k+1)}
=\frac{\pi(2k-1)!!}{2^{k+1}k!}\gamma^{-(2k+1)}.
\]
Substitution proves \eqref{eq:mainoperatorintegral}.
\end{proof}

The first three odd instances, in the same quick-reference form as the even examples in Section~\ref{sec:evenintegral}, are
\begin{align}
Z(3)
&=\frac{1}{\pi}\int_{0}^{\infty}
\left(
\frac{L'(x)}{x^3}
-\frac{L''(x)}{x^2}
\right)\dd x,
\label{eq:Z3quick}\\
Z(5)
&=\frac{1}{3\pi}\int_{0}^{\infty}
\left(
\frac{L'''(x)}{x^3}
-\frac{3L''(x)}{x^4}
+\frac{3L'(x)}{x^5}
\right)\dd x,
\label{eq:Z5quick}\\
Z(7)
&=-\frac{1}{15\pi}\int_{0}^{\infty}
\left(
\frac{L^{(4)}(x)}{x^4}
-\frac{6L'''(x)}{x^5}
+\frac{15L''(x)}{x^6}
-\frac{15L'(x)}{x^7}
\right)\dd x.
\label{eq:Z7quick}
\end{align}
Each parenthesized combination must be kept intact, since its apparent
singular terms cancel at $x=0$.

\section{Reference values \texorpdfstring{$Z(2)$}{Z(2)} to \texorpdfstring{$Z(7)$}{Z(7)}}
\label{sec:referencevalues}

This section collects the first six cases in one place.  The formulas are
intended as a convenient reference for direct numerical evaluation and for
checking symbolic implementations of the operators.

\subsection{\texorpdfstring{$\Dop_nL$}{D n L} and \texorpdfstring{$\xi$}{xi}}
\label{subsec:DnLxi}

The first three operator expansions are
\begin{align}
\Dop_1L
&=x^2L''-xL',
\label{eq:D1expanded}\\
\Dop_2L
&=x^3L'''-3x^2L''+3xL',
\label{eq:D2expanded}\\
\Dop_3L
&=x^4L^{(4)}-6x^3L'''+15x^2L''-15xL'.
\label{eq:D3expanded}
\end{align}
All derivatives in these formulas are evaluated at the integration variable
$x$.  Primes below denote differentiation of the Riemann $\xi$-function
with respect to its argument.  Since the denominator $\xi(\tfrac12)$ in the
definition of $L$ is constant,
\begin{align}
L'(x)
&=\frac{\xi'(x+\tfrac12)}{\xi(x+\tfrac12)},
\label{eq:LprimeXi}\\
L''(x)
&=\frac{\xi''(x+\tfrac12)}{\xi(x+\tfrac12)}
-\frac{\bigl[\xi'(x+\tfrac12)\bigr]^2}
{\bigl[\xi(x+\tfrac12)\bigr]^2},
\label{eq:LsecondXi}\\
L'''(x)
&=\frac{\xi'''(x+\tfrac12)}{\xi(x+\tfrac12)}
-\frac{3\,\xi'(x+\tfrac12)\xi''(x+\tfrac12)}
{\bigl[\xi(x+\tfrac12)\bigr]^2}
+\frac{2\bigl[\xi'(x+\tfrac12)\bigr]^3}
{\bigl[\xi(x+\tfrac12)\bigr]^3},
\label{eq:LthirdXi}\\
L^{(4)}(x)
&=\frac{\xi^{(4)}(x+\tfrac12)}{\xi(x+\tfrac12)}
-\frac{
4\,\xi'(x+\tfrac12)\xi'''(x+\tfrac12)
+3\bigl[\xi''(x+\tfrac12)\bigr]^2
}{\bigl[\xi(x+\tfrac12)\bigr]^2}
\notag\\
&\quad
+\frac{12\bigl[\xi'(x+\tfrac12)\bigr]^2\xi''(x+\tfrac12)}
{\bigl[\xi(x+\tfrac12)\bigr]^3}
-\frac{6\bigl[\xi'(x+\tfrac12)\bigr]^4}
{\bigl[\xi(x+\tfrac12)\bigr]^4}.
\label{eq:LfourthXi}
\end{align}

\subsection{The value \texorpdfstring{$Z(2)$}{Z(2)}}
\label{subsec:referenceZ2}

The centered Taylor-coefficient formula gives
\begin{equation}
\boxed{
Z(2)=\frac12L''(0)
=\frac12\frac{\xi''(\tfrac12)}{\xi(\tfrac12)},
}
\label{eq:Z2localreference}
\end{equation}
where $\xi'(\tfrac12)=0$ follows from the functional equation.  The
convergent integral form is
\begin{equation}
\boxed{
Z(2)=\frac12\int_0^\infty
\left(
\frac{L'(x)}{x^2}-\frac{L''(x)}{x}
\right)\dd x.
}
\label{eq:Z2referenceintegral}
\end{equation}
Using \eqref{eq:LprimeXi} and \eqref{eq:LsecondXi},
\begin{equation}
\boxed{
Z(2)=\frac12\int_0^\infty\Bigg(
\frac{\xi'(x+\tfrac12)}{x^2\xi(x+\tfrac12)}
-\frac{\xi''(x+\tfrac12)}{x\xi(x+\tfrac12)}
+\frac{\bigl[\xi'(x+\tfrac12)\bigr]^2}
{x\bigl[\xi(x+\tfrac12)\bigr]^2}
\Bigg)\dd x.
}
\label{eq:Z2ordinaryxi}
\end{equation}

\subsection{The value \texorpdfstring{$Z(3)$}{Z(3)}}
\label{subsec:referenceZ3}

For $k=1$, \eqref{eq:remainderintegral} becomes
\begin{equation}
\boxed{
Z(3)=\frac{3}{\pi}\int_{0}^{\infty}
\frac{Z(2)x^2-L(x)}{x^4}\,\dd x.
}
\label{eq:Z3remainder}
\end{equation}
Using only the defining $\xi$-function ratio for $L(x)$, this is
\begin{equation}
\boxed{
Z(3)=\frac{3}{\pi}\int_{0}^{\infty}
\frac{
Z(2)x^2-
\log\!\left(\dfrac{\xi(x+\tfrac12)}{\xi(\tfrac12)}\right)
}{x^4}\,\dd x.
}
\label{eq:Z3xi}
\end{equation}
The expanded operator formula is
\begin{equation}
\boxed{
Z(3)=\frac{1}{\pi}\int_{0}^{\infty}
\left(
\frac{L'(x)}{x^3}
-\frac{L''(x)}{x^2}
\right)\dd x.
}
\label{eq:Z3derivative}
\end{equation}
The two terms in parentheses are to be combined before integration; their
singular parts at $x=0$ cancel exactly.  Using
\eqref{eq:LprimeXi}--\eqref{eq:LsecondXi} gives
\begin{equation}
\boxed{
Z(3)=\frac{1}{\pi}\int_{0}^{\infty}\Bigg(
\frac{\xi'(x+\tfrac12)}{x^3\,\xi(x+\tfrac12)}
-\frac{\xi''(x+\tfrac12)}{x^2\,\xi(x+\tfrac12)}
+\frac{\bigl[\xi'(x+\tfrac12)\bigr]^2}
{x^2\bigl[\xi(x+\tfrac12)\bigr]^2}
\Bigg)\dd x.
}
\label{eq:Z3ordinaryxi}
\end{equation}

\subsection{The value \texorpdfstring{$Z(4)$}{Z(4)}}
\label{subsec:referenceZ4}

At the center, the functional equation makes all odd derivatives of $\xi$
vanish.  Hence
\begin{equation}
\boxed{
Z(4)=-\frac1{12}L^{(4)}(0)
=-\frac1{12}\left[
\frac{\xi^{(4)}(\tfrac12)}{\xi(\tfrac12)}
-3\left(\frac{\xi''(\tfrac12)}{\xi(\tfrac12)}\right)^2
\right].
}
\label{eq:Z4localreference}
\end{equation}
The convergent integral form is
\begin{equation}
\boxed{
Z(4)=\frac14\int_0^\infty
\left(
\frac{L'''(x)}{x^2}
-\frac{3L''(x)}{x^3}
+\frac{3L'(x)}{x^4}
\right)\dd x.
}
\label{eq:Z4referenceintegral}
\end{equation}
Substitution of \eqref{eq:LprimeXi}--\eqref{eq:LthirdXi} yields
\begin{equation}
\boxed{
\begin{aligned}
Z(4)=\frac14\int_0^\infty\Bigg(&
\frac{\xi'''(x+\tfrac12)}{x^2\xi(x+\tfrac12)}
-\frac{3\xi'(x+\tfrac12)\xi''(x+\tfrac12)}
{x^2\bigl[\xi(x+\tfrac12)\bigr]^2}\\
&+\frac{2\bigl[\xi'(x+\tfrac12)\bigr]^3}
{x^2\bigl[\xi(x+\tfrac12)\bigr]^3}
-\frac{3\xi''(x+\tfrac12)}{x^3\xi(x+\tfrac12)}\\
&+\frac{3\bigl[\xi'(x+\tfrac12)\bigr]^2}
{x^3\bigl[\xi(x+\tfrac12)\bigr]^2}
+\frac{3\xi'(x+\tfrac12)}{x^4\xi(x+\tfrac12)}
\Bigg)\dd x.
\end{aligned}
}
\label{eq:Z4ordinaryxi}
\end{equation}

\subsection{The value \texorpdfstring{$Z(5)$}{Z(5)}}
\label{subsec:referenceZ5}

For $k=2$,
\begin{equation}
\boxed{
Z(5)=\frac{5}{\pi}\int_{0}^{\infty}
\frac{
L(x)-Z(2)x^2+\tfrac12Z(4)x^4
}{x^6}\,\dd x.
}
\label{eq:Z5remainder}
\end{equation}
Equivalently, entirely in terms of the $\xi$-function,
\begin{equation}
\boxed{
Z(5)=\frac{5}{\pi}\int_{0}^{\infty}
\frac{
\log\!\left(\dfrac{\xi(x+\tfrac12)}{\xi(\tfrac12)}\right)
-Z(2)x^2+\tfrac12Z(4)x^4
}{x^6}\,\dd x.
}
\label{eq:Z5xi}
\end{equation}
The expanded operator formula is
\begin{equation}
\boxed{
Z(5)=\frac{1}{3\pi}\int_{0}^{\infty}
\left(
\frac{L'''(x)}{x^3}
-\frac{3L''(x)}{x^4}
+\frac{3L'(x)}{x^5}
\right)\dd x.
}
\label{eq:Z5derivative}
\end{equation}
Again, the displayed combination is integrable even though its individual
summands are not separately integrable at $x=0$.  Using
\eqref{eq:LprimeXi}--\eqref{eq:LthirdXi}, this becomes
\begin{equation}
\boxed{
\begin{aligned}
Z(5)=\frac{1}{3\pi}\int_{0}^{\infty}\Bigg(&
\frac{\xi'''(x+\tfrac12)}{x^3\,\xi(x+\tfrac12)}
-\frac{3\,\xi'(x+\tfrac12)\xi''(x+\tfrac12)}
{x^3\bigl[\xi(x+\tfrac12)\bigr]^2}\\
&+\frac{2\bigl[\xi'(x+\tfrac12)\bigr]^3}
{x^3\bigl[\xi(x+\tfrac12)\bigr]^3}
-\frac{3\,\xi''(x+\tfrac12)}{x^4\,\xi(x+\tfrac12)}\\
&+\frac{3\bigl[\xi'(x+\tfrac12)\bigr]^2}
{x^4\bigl[\xi(x+\tfrac12)\bigr]^2}
+\frac{3\,\xi'(x+\tfrac12)}{x^5\,\xi(x+\tfrac12)}
\Bigg)\dd x.
\end{aligned}
}
\label{eq:Z5ordinaryxi}
\end{equation}

\subsection{The value \texorpdfstring{$Z(6)$}{Z(6)}}
\label{subsec:referenceZ6}

The local centered formula may be written
\begin{equation}
\boxed{
Z(6)=\frac1{240}L^{(6)}(0)
=\frac1{240}\Bigg[
\frac{\xi^{(6)}(\tfrac12)}{\xi(\tfrac12)}
-\frac{15\xi''(\tfrac12)\xi^{(4)}(\tfrac12)}
{\bigl[\xi(\tfrac12)\bigr]^2}
+\frac{30\bigl[\xi''(\tfrac12)\bigr]^3}
{\bigl[\xi(\tfrac12)\bigr]^3}
\Bigg].
}
\label{eq:Z6localreference}
\end{equation}
The convergent integral form is
\begin{equation}
\boxed{
Z(6)=-\frac1{16}\int_0^\infty
\left(
\frac{L^{(4)}(x)}{x^3}
-\frac{6L'''(x)}{x^4}
+\frac{15L''(x)}{x^5}
-\frac{15L'(x)}{x^6}
\right)\dd x.
}
\label{eq:Z6referenceintegral}
\end{equation}
Using \eqref{eq:LprimeXi}--\eqref{eq:LfourthXi},
\begin{equation}
\boxed{
\begin{aligned}
Z(6)=-\frac1{16}\int_0^\infty\Bigg(&
\frac{\xi^{(4)}(x+\tfrac12)}{x^3\xi(x+\tfrac12)}
-\frac{4\xi'(x+\tfrac12)\xi'''(x+\tfrac12)}
{x^3\bigl[\xi(x+\tfrac12)\bigr]^2}\\
&-\frac{3\bigl[\xi''(x+\tfrac12)\bigr]^2}
{x^3\bigl[\xi(x+\tfrac12)\bigr]^2}
+\frac{12\bigl[\xi'(x+\tfrac12)\bigr]^2\xi''(x+\tfrac12)}
{x^3\bigl[\xi(x+\tfrac12)\bigr]^3}\\
&-\frac{6\bigl[\xi'(x+\tfrac12)\bigr]^4}
{x^3\bigl[\xi(x+\tfrac12)\bigr]^4}
-\frac{6\xi'''(x+\tfrac12)}{x^4\xi(x+\tfrac12)}\\
&+\frac{18\xi'(x+\tfrac12)\xi''(x+\tfrac12)}
{x^4\bigl[\xi(x+\tfrac12)\bigr]^2}
-\frac{12\bigl[\xi'(x+\tfrac12)\bigr]^3}
{x^4\bigl[\xi(x+\tfrac12)\bigr]^3}\\
&+\frac{15\xi''(x+\tfrac12)}{x^5\xi(x+\tfrac12)}
-\frac{15\bigl[\xi'(x+\tfrac12)\bigr]^2}
{x^5\bigl[\xi(x+\tfrac12)\bigr]^2}\\
&-\frac{15\xi'(x+\tfrac12)}{x^6\xi(x+\tfrac12)}
\Bigg)\dd x.
\end{aligned}
}
\label{eq:Z6ordinaryxi}
\end{equation}

\subsection{The value \texorpdfstring{$Z(7)$}{Z(7)}}
\label{subsec:referenceZ7}

For $k=3$,
\begin{equation}
\boxed{
Z(7)=\frac{7}{\pi}\int_{0}^{\infty}
\frac{
Z(2)x^2-\tfrac12Z(4)x^4
+\tfrac13Z(6)x^6-L(x)
}{x^8}\,\dd x.
}
\label{eq:Z7remainder}
\end{equation}
Thus the direct $\xi$-function form is
\begin{equation}
\boxed{
Z(7)=\frac{7}{\pi}\int_{0}^{\infty}
\frac{
Z(2)x^2-\tfrac12Z(4)x^4
+\tfrac13Z(6)x^6
-\log\!\left(\dfrac{\xi(x+\tfrac12)}{\xi(\tfrac12)}\right)
}{x^8}\,\dd x.
}
\label{eq:Z7xi}
\end{equation}
The expanded operator formula is
\begin{equation}
\boxed{
Z(7)=-\frac{1}{15\pi}\int_{0}^{\infty}
\left(
\frac{L^{(4)}(x)}{x^4}
-\frac{6L'''(x)}{x^5}
+\frac{15L''(x)}{x^6}
-\frac{15L'(x)}{x^7}
\right)\dd x.
}
\label{eq:Z7derivative}
\end{equation}
As before, the cancellations must be retained at the level of the full
parenthesized expression.  Using
\eqref{eq:LprimeXi}--\eqref{eq:LfourthXi} gives
\begin{equation}
\boxed{
\begin{aligned}
Z(7)=-\frac{1}{15\pi}\int_{0}^{\infty}\Bigg(&
\frac{\xi^{(4)}(x+\tfrac12)}{x^4\,\xi(x+\tfrac12)}
-\frac{4\,\xi'(x+\tfrac12)\xi'''(x+\tfrac12)}
{x^4\bigl[\xi(x+\tfrac12)\bigr]^2}\\
&-\frac{3\bigl[\xi''(x+\tfrac12)\bigr]^2}
{x^4\bigl[\xi(x+\tfrac12)\bigr]^2}
+\frac{12\bigl[\xi'(x+\tfrac12)\bigr]^2\xi''(x+\tfrac12)}
{x^4\bigl[\xi(x+\tfrac12)\bigr]^3}\\
&-\frac{6\bigl[\xi'(x+\tfrac12)\bigr]^4}
{x^4\bigl[\xi(x+\tfrac12)\bigr]^4}
-\frac{6\,\xi'''(x+\tfrac12)}{x^5\,\xi(x+\tfrac12)}\\
&+\frac{18\,\xi'(x+\tfrac12)\xi''(x+\tfrac12)}
{x^5\bigl[\xi(x+\tfrac12)\bigr]^2}
-\frac{12\bigl[\xi'(x+\tfrac12)\bigr]^3}
{x^5\bigl[\xi(x+\tfrac12)\bigr]^3}\\
&+\frac{15\,\xi''(x+\tfrac12)}{x^6\,\xi(x+\tfrac12)}
-\frac{15\bigl[\xi'(x+\tfrac12)\bigr]^2}
{x^6\bigl[\xi(x+\tfrac12)\bigr]^2}\\
&-\frac{15\,\xi'(x+\tfrac12)}{x^7\,\xi(x+\tfrac12)}
\Bigg)\dd x.
\end{aligned}
}
\label{eq:Z7ordinaryxi}
\end{equation}

\section{Recovering the first zero ordinate \texorpdfstring{$\gamma_1$}{gamma1}}
\label{sec:gamma1}

Let
\[
0<\gamma_1<\gamma_2<\gamma_3<\cdots
\]
denote the distinct positive ordinates, and let $m_n$ be the multiplicity of
$\gamma_n$.  For real $s>1$, factor the contribution of the first ordinate
from the defining series:
\begin{equation}
Z(s)=m_1\gamma_1^{-s}\bigl(1+R_s\bigr),
\qquad
R_s=\sum_{n\geq2}\frac{m_n}{m_1}
\left(\frac{\gamma_1}{\gamma_n}\right)^s.
\label{eq:Zdominatedfirstzero}
\end{equation}
Since $R_s\to0$ as $s\to+\infty$,
\begin{equation}
\frac{Z(s)}{Z(s+1)}
=\gamma_1\frac{1+R_s}{1+R_{s+1}}
\longrightarrow\gamma_1.
\label{eq:gamma1generalratio}
\end{equation}
In particular, the even and odd values developed above give
\begin{equation}
\boxed{
\gamma_1=\lim_{k\to\infty}\frac{Z(2k)}{Z(2k+1)}.
}
\label{eq:gamma1evenoddlimit}
\end{equation}
Substituting the even- and odd-value operator integrals into
\eqref{eq:gamma1evenoddlimit} gives, with the same operator
$\Dop_k=\Eop\prod_{j=1}^{k}(\Eop-2j)$,
\begin{equation}
\boxed{
Z(2k)=\frac{(-1)^k}{2^k(k-1)!}
\int_{0}^{\infty}
\frac{\Dop_kL(x)}{x^{2k+1}}\,\dd x,
}
\label{eq:gamma1evenintegral}
\end{equation}
and
\begin{equation}
\boxed{
Z(2k+1)=\frac{(-1)^k}{\pi(2k-1)!!}
\int_{0}^{\infty}
\frac{\Dop_kL(x)}{x^{2k+2}}\,\dd x.
}
\label{eq:gamma1oddintegral}
\end{equation}
The factors $(-1)^k$ cancel in the quotient, while
\[
\frac{\pi(2k-1)!!}{2^k(k-1)!}
=\frac{\pi k}{4^k}\binom{2k}{k}
=\sqrt{\pi}\,\frac{\Gamma(k+\tfrac12)}{\Gamma(k)}.
\]
Therefore the first ordinate has the direct integral-ratio representation
\begin{equation}
\boxed{
\gamma_1
=\lim_{k\to\infty}\:
\sqrt{\pi}\,\frac{\Gamma(k+\tfrac12)}{\Gamma(k)}\:
\frac{
\displaystyle\int_{0}^{\infty}
\dfrac{\Dop_kL(x)}{x^{2k+1}}\,\dd x
}{
\displaystyle\int_{0}^{\infty}
\dfrac{\Dop_kL(x)}{x^{2k+2}}\,\dd x
}.
}
\label{eq:gamma1integralratio}
\end{equation}
Moreover, $R_s>R_{s+1}$, so every finite ratio in
\eqref{eq:gamma1generalratio} lies above $\gamma_1$.

In Table~\ref{tab:gamma1ratios} we compute the ratios to high-precision using the Arias de Reyna algorithm, see \cite{Arias2020}\cite{Kawalec2024}, as numerical computation of these integrals directly to high precision is difficult.

\begin{table}[H]
\centering
\caption{Convergence of $Z(2k)/Z(2k+1)$ to the first positive ordinate.
The third column counts the consecutive digits after the decimal point
that agree with $\gamma_1$; the last agreeing decimal digit is
highlighted in blue.}
\label{tab:gamma1ratios}
\begingroup
\setlength{\tabcolsep}{12pt}
\fontsize{12pt}{12pt}\selectfont
\renewcommand{\arraystretch}{1.12}
\begin{tabular}{@{}c>{\ttfamily}l c@{}}
\toprule
$k$
& \multicolumn{1}{c}{$Z(2k)/Z(2k+1)$}
& Matching decimal digits \\
\midrule
1
& 31.67027320838126
& 0 \\

2
& 16.660446197833857
& 0 \\

3
& 14.901167680867425
& 0 \\

4
& 14.416266812062819
& 0 \\

5
& 14.247171591098658
& 0 \\

6
& 14.\textcolor{blue}{1}81687529802197
& 1 \\

7
& 14.\textcolor{blue}{1}54865750419131
& 1 \\

8
& 14.\textcolor{blue}{1}43509103505405
& 1 \\

9
& 14.1\textcolor{blue}{3}8598821034345
& 2 \\

10
& 14.1\textcolor{blue}{3}6446297126356
& 2 \\

25
& 14.134725\textcolor{blue}{1}52861886
& 7 \\

50
& 14.13472514173469\textcolor{blue}{3}817182033
& 15 \\

100
& 14.13472514173469379045725198356247\textcolor{blue}{0}4250132
& 33 \\
\bottomrule
\end{tabular}
\endgroup
\end{table}

The limiting value is
\[
\gamma_1=14.134725141734693790457251983562470270\ldots,
\]
so the ratios obtained from the tabulated values already show the expected
convergence from above.

\section{Conclusion}

The single logarithm
\[
L(x)=\log\!\left(\frac{\xi(x+\tfrac12)}{\xi(\tfrac12)}\right)
\]
provides a common source for both parities of the secondary-zeta values.  Its Taylor coefficients at $x=0$ give $Z(2k)$ directly, while its regularized Mellin remainders give $Z(2k+1)$.  The centered Euler operator $\Dop_k$ packages the necessary Taylor subtraction into one differential expression and converts each logarithmic zero factor into a positive resolvent kernel.  The reference formulas for $Z(2)$ through $Z(7)$ may therefore be written either with derivatives of $L$ or entirely with ordinary derivatives $\xi^{(j)}(x+\tfrac12)$, with every shifted argument displayed explicitly.  Finally, the high-order values recover the first ordinate through $\gamma_1=\lim_{k\to\infty}Z(2k)/Z(2k+1)$.

\section*{AI Disclosure}
This manuscript was developed with GPT~5.6 Pro (OpenAI), with additional assistance
from Claude Opus~5.5 (Anthropic) and Grok~4.7 (SpaceXAI). These systems assisted with
drafting the text, deriving equations, generating and interpreting numerical
simulations, and identifying candidate literature. The core result of this paper,
the differential operator integral representation, was found by GPT~5.6 Pro. The
author directed the research, edited the final manuscript, and independently
reviewed and verified all mathematical arguments, numerical computations, and
references. The author takes full responsibility for the content and correctness of
the paper. No AI system is listed as an author.

\end{document}